# About some generalizations trigonometric splines

Denysiuk V.P.
Dr. of Ps.-M.. Sciences, Professor, Kiev, Ukraine
National Aviation University
kvomden@nau.edu.ua

## Annotation

Methods of constructing trigonometric fundamental splines with constant sign and sign-changing convergence factors are given. An example and graphics illustrating the concepts of convergence and interpolation grids are given. Some methods of constructing constant-sign and sign-changing coefficients of convergence of trigonometric splines are considered.

## Introduction

Approximation, respectively representation, of an arbitrary known or unknown function through a set of some special functions can be considered as a central topic of analysis. We will use the term "special functions" to refer to classes of algebraic and trigonometric polynomials and their modifications; at the same time, we believe that the classes of trigonometric polynomials also include trigonometric series. As a rule, such special functions are easy to calculate and have interesting analytical properties [1].

One of the most successful modifications of algebraic polynomials are polynomial splines that are stitched together from segments of these polynomials according to a certain scheme. The theory of polynomial splines appeared relatively recently and is well developed (see, for example, [2], [3], [4] [5], etc.). The advantages of polynomial splines include the fact that they can be given certain smoothness properties, as well as their approximate properties [6]. The main disadvantage of polynomial splines, in our opinion, is their piecemeal structure, which greatly complicates their use in analytical transformations.

Later it turned out [7], [8] that there are also modifications of trigonometric series whose sums depend on several parameters and have the same properties as polynomial splines [9]; moreover, the class of such modified series is quite broad and includes the class of polynomial periodic splines. This gave reason to call the class of such series trigonometric interpolation splines.

Convergence of trigonometric series that determine trigonometric interpolation splines, provided by convergence factors [7], [8], which have the order of decreasing $O\left(k^{-(1+r)}\right)$, ($r > 0$); since these series coincide uniformly, they are trigonometric Fourier series with special coefficients. In this work, we will limit ourselves to consideration of integer values of the parameter $r$ (i.e. case $r = 1, 2, ...$); note that trigonometric splines of fractional powers (that is, of non-integer parameter values $r$) was considered in [10].

Trigonometric series that provide trigonometric interpolation splines can naturally be divided into two components - even and odd. In turn, each of these components can be broken down into three more components: low-frequency, medium-frequency, and high-frequency. These components can be considered together with the coefficients $\Gamma = \{\gamma_1, \gamma_2, \gamma_3\}$ and $H = \{\eta_1, \eta_2, \eta_3\}$, where the components of the vectors $\Gamma$ and $H$ real numbers and at least one of the components $\gamma_2, \gamma_3, \eta_2, \eta_3$ not equal to 0. It is clear that with vectors

$\Gamma = \{\gamma_1, 0, 0\}$ and $H = \{\eta_1, 0, 0\}$ we have a trigonometric interpolation polynomial. Note that with the help of vectors $\Gamma$ and $H$ it is easy to model low-, medium-, and high-pass filters, which is important in the problems of the theory of modeling and signal processing.

Trigonometric splines with vectors $\Gamma = \{1,1,1\}$ and $H = \{1,1,1\}$ we will call simple. In the notation of simple trigonometric splines, we will omit their dependence on vectors $\Gamma$ and $H$.

Trigonometric interpolation splines assume a number of generalizations; below, we consider only generalizations that relate to types of convergence factors.

**The goal of the work.**

1. Consideration of methods for constructing one- and two-parameter convergence factors of trigonometric fundamental simple splines.

**Main part.**

**General methods of constructing trigonometric fundamental splines with constant coefficients of convergence.**

Let $W_v^r$ - class $2\pi$-periodic functions having a completely continuous derivative of order $r-1$ ($r = 1, 2, ...$), and the derivative is of order $r$ is a function of bounded variation; it is clear that $W_v^1$ is a class of continuous functions. A symbol $W_v^0$ denote the class of piecewise constant functions with a finite number of discontinuity points.

At a stretch $[0, 2\pi)$ we will consider grids $\Delta_N^{(I)} = \{t_j^{(I)}\}_{j=1}^N$, ($I = 0,1$), $t_j^{(0)} = \frac{2\pi}{N}(j-1)$, $t_j^{(1)} = \frac{\pi}{N}(j-1)$. Let also on $[0, 2\pi)$ given a periodic function $f(t) \in W_v^r$, and its values are known $f_j^{(I)} = \{f(t_j^{(I)})\}_{j=1}^N$ in grid nodes $\Delta_N^{(I)}$.

Let's introduce the concepts of stitching grids and interpolation of trigonometric splines.

We will call the mesh a stitching mesh $\Delta_N^{(I)}$, on which the polynomial analog of the trigonometric spline of the 0th degree has discontinuities of the 1st kind of the jump type or the polynomial analogs of the trigonometric splines of the 1st degree are stitched together; indicator $I$ we will mark the stitching grids through $I_1$, ($I_1 = 0,1$). Note that trigonometric splines of degree 0,1 were chosen for the reasons that in these cases the points of discontinuity or splicing appear in the most relief.

We will call the grid an interpolation grid or an interpolation grid $\Delta_N^{(I)}$, in the nodes of which interpolation of the approximated function is carried out by trigonometric splines; indicator $I$ interpolation grids will be denoted by $I_2$, ($I_2 = 0,1$).

As we said before [11], the system $st_k(r,t) \in W_v^{r-1}$, $r = 1, 2, ..., N$ trigonometric fundamental splines on an interpolation grid $\Delta_N^{(I)}$ meets the conditions: $st_k^{(I)}(r, t_j^{(I)}) = \begin{cases} 1, k = j; \\ 0, k \neq j. \end{cases}$ ($k, j = 1, ..., N$). (1)

Using a system of trigonometric fundamental splines $st_k^{(I)}(r,t), k = 1, 2, \cdots, N$, trigonometric interpolation spline can be written as: $St_n^{(I)}(r,t) = \sum_{k=1}^N f_k^{(I)} st_k^{(I)}(r,t)$. (2)

Note that in the definition of trigonometric fundamental splines we are talking about interpolation grids. Next, stitch and interpolation grid indices, as well as vectors $\Gamma$ and $H$ we will enter splines in notation.

Parameter $r$, ($r = 1, 2, ...$), determines the order of decreasing convergence factors $\sigma(r,j)$. The most natural as convergence multipliers are the constant multipliers: $\sigma 0(r,j) = \left(\frac{1}{j}\right)^{1+r}$, ($j = 1, 2, ...$). (3)

It is clear that the parameter $r$, ($r = 1, 2, ...$) determines the number of continuous derivatives of trigonometric fundamental splines.

Note that otherwise, we considered trigonometric splines with sign-changing convergence factors of the Riemann type [9]; in this work, we consider constant coefficients of convergence (3), the only property of which is the order of descent. It is clear that sign-changing and sign-constant convergence factors lead to different algorithms for constructing trigonometric splines, which we will present below.

Functions of the system of fundamental trigonometric splines $\{ts_k(I_1,I_2,\Gamma,H,\sigma,r,t)\}_{k=1}^{N}$ with stitching grid $I_1$, interpolation grid $I_1$, vectors $\Gamma$ and $H$, belonging to space $W_\nu^{r-1}$, with sign constant convergence factors $\sigma 0(r,j)$ can be submitted in the form of

$$st_k(I_1,I_2,\Gamma,H,\sigma 0,r,N,t) =$$

$$= \frac{1}{N}\left(1 + 2\sum_{j=1}^{\frac{N-1}{2}}\left(\frac{c_j(I_1,I_2,\Gamma,\sigma 0,r,k,N,t)}{hc_j(I_1,I_2,\Gamma,\sigma 0,r,k,N)} + \frac{s_j(I_1,I_2,H,\sigma 0,r,k,N,t)}{hs_j(I_1,I_2,H,\sigma 0,r,k,N)}\right)\right), \quad (4)$$

where

$$c_j(I_1,I_2,\Gamma,\sigma 0,r,k,t) = \gamma_1\sigma 0(r,j)\cos(jt)\cos(j(x_k^{(0)})^{1-I_2}(x_k^{(1)})^{I_2}) +$$

$$+\sum_{m=1}^{\infty}(-1)^{m(I_1+I_2)}\left[\gamma_2\sigma 0(r,mN-j)(-1)^{1+r}\cos((mN-j)t)\cos((mN-j)(x_k^{(0)})^{1-I_2}(x_k^{(1)})^{I_2}) + (5)\right.$$

$$\left.+\gamma_3\sigma 0(r,mN+j)\cos((mN+j)t)\cos((mN+j)(x_k^{(0)})^{1-I_2}(x_k^{(1)})^{I_2})\right];$$

$$hc_j(I_1,I_2,\Gamma,\sigma 0,r,N) = \gamma_1\sigma 0(r,j) + \sum_{m=1}^{\infty}(-1)^{m(I_1+I_2)}\left[\gamma_2(-1)^{1+r}\sigma 0(r,mN-j) + \gamma_3\sigma 0(r,mN+j)\right]; \quad (6)$$

$$s_j(I_1,I_2,H,\sigma 0,r,k,t) = \eta_1\sigma 0(r,j)\sin(jt)\sin(j(x_k^{(0)})^{1-I_2}(x_k^{(1)})^{I_2}) +$$

$$+\sum_{m=1}^{\infty}(-1)^{m(I_1+I_2)}\left[\eta_2\sigma 0(r,mN-j)(-1)^{1+r}\sin((mN-j)t)\sin((mN-j)(x_k^{(0)})^{1-I_2}(x_k^{(1)})^{I_2}) + (7)\right.$$

$$\left.+\eta_3\sigma 0(r,mN+j)\sin((mN+j)t)\sin((mN+j)(x_k^{(0)})^{1-I_2}(x_k^{(1)})^{I_2})\right];$$

$$hs_j(I_1,I_2,H,\sigma 0,r,N) = \eta_1\sigma 0(r,j) + \sum_{m=1}^{\infty}(-1)^{m(I_1+I_2)}\left[\eta_2(-1)^{1+r}\sigma 0(r,mN-j) + \eta_3\sigma 0(r,mN+j)\right]. \quad (8)$$

The construction of the function of the system of trigonometric fundamental simple splines is carried out according to the formulas: $st_k(I_1,I_2,\sigma 0,r,N,t) = \dfrac{1}{N}\left(1 + 2\sum_{j=1}^{\frac{N-1}{2}}\dfrac{c_j(I_1,I_2,\sigma 0,r,k,N,t)}{h_j(I_1,I_2,\sigma 0,r,N)}\right), \quad (9)$

where

$$c_j(I_1,I_2,\sigma 0,r,k,t) = \sigma 0(r,j)\cos(j(t-(x_k^{(0)})^{1-I_2}(x_k^{(1)})^{I_2})) +$$

$$+\sum_{m=1}^{\infty}(-1)^{m(I_1+I_2)}\left[\sigma 0(r,mN-j)(-1)^{1+r}\cos((mN-j)(t-(x_k^{(0)})^{1-I_2}(x_k^{(1)})^{I_2})) + (10)\right.$$

$$\left.+\sigma 0(r,mN+j)\cos((mN+j)(t-(x_k^{(0)})^{1-I_2}(x_k^{(1)})^{I_2}))\right];$$

$$h_j(I_1,I_2,\Gamma,\sigma 0,r,N) = \sigma 0(r,j) + \sum_{m=1}^{\infty}(-1)^{m(I_1+I_2)}\left[(-1)^{1+r}\sigma 0(r,mN-j) + \sigma 0(r,mN+j)\right]. \quad (11)$$

In this work, we will limit ourselves to considering the system of trigonometric fundamental simple splines only.

To illustrate the introduced concepts we will give an example of stitching grids and interpolation grids.

**Example.** Let's put $N=9$, $r=1$, $k=5$. Let us present the graphs of the trigonometric fundamental simple spline $st_5(I_1,I_2,\sigma 0,1,9,t)$ for different meshes of stitching and interpolation. Note that in the graphs below, the vertical grid lines are drawn at the grid nodes $\Delta_9^{(0)}$; respectively, grid nodes $\Delta_9^{(1)}$ located between grid nodes $\Delta_9^{(0)}$.

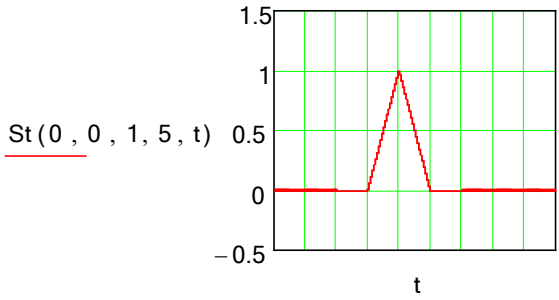

Fig. 1. Graph of trigonometric fundamental simple spline $st_5(0,0,\sigma 0,1,9,t)$; stitching and interpolation is carried out at grid nodes $\Delta_9^{(0)}$.

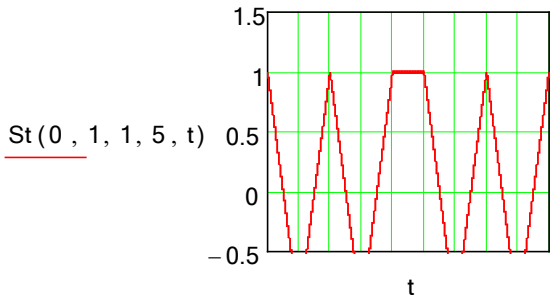

Fig. 2. Graph of trigonometric fundamental simple spline $st_5(0,1,\sigma 0,1,9,t)$; stitching in grid nodes $\Delta_9^{(0)}$, interpolation at grid nodes $\Delta_9^{(1)}$.

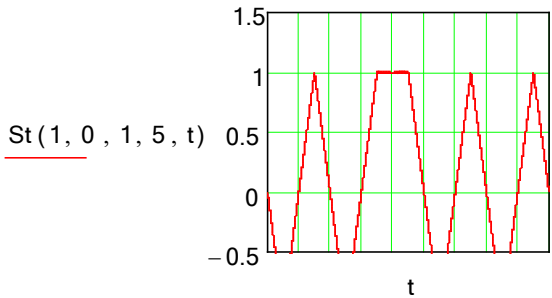

Fig. 3. Graph of trigonometric fundamental simple spline $st_5(1,0,\sigma 0,1,9,t)$; stitching in grid nodes $\Delta_9^{(1)}$, interpolation at grid nodes $\Delta_9^{(0)}$.

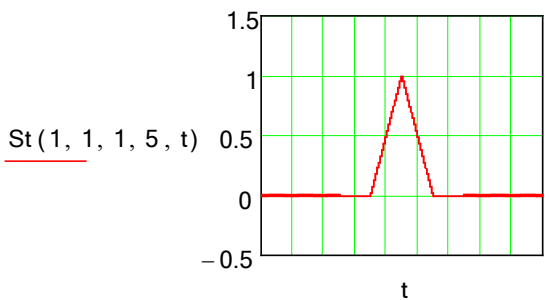

Fig. 4. Graph of trigonometric fundamental simple spline $st_5(1,1,\sigma 0,1,9,t)$; stitching and interpolation at grid nodes $\Delta_9^{(1)}$.

Note that only for splines $st_5(0,0,\sigma 0,2r-1,9,t)$ and $st_5(0,1,\sigma 0,2r,9,t)$ ($r=1,2,...$ known polynomial analogues; for other splines, polynomial analogues are unknown, but they are easy to construct at least at $r=1$.

Let's consider some methods of constructing factors of convergence of trigonometric splines.

1. **The method of construction of constant coefficients of convergence, , is based on the equivalence of infinitesimals.**

A more general type of constant coefficients $\sigma 0(r,j)$ are convergence factors $\sigma 0(\alpha,r,j)$, which are defined as follows $\sigma 0(\alpha,r,j) = \dfrac{\alpha}{j^{1+r}}$, ($j = 1,2,...$), (12) where $\alpha$, ($\alpha > 0$) is some parameter.

Note that the coefficients of convergence are constant $\sigma 0(\alpha,r,j)$ at $j \to \infty$ are infinitely small quantities. Using the known equivalence relations of infinitesimals, it is easy to construct the following convergence factors:

$$\sigma 01(\alpha,r,j) = \sin(\sigma 0(\alpha,r,j));$$

$$\sigma 02(\alpha,r,j) = \arcsin(\sigma 0(\alpha,r,j));$$

$$\sigma 03(\alpha,r,j) = tg(\sigma 0(\alpha,r,j))$$

$$\sigma 04(\alpha,r,j) = arctg(\sigma 0(\alpha,r,j))$$

$$\sigma 5(\alpha,a,r,j) = \ln a \log_a(1+\sigma 0(\alpha,r,j))$$

$$\sigma 06(\alpha,a,r,j) = \dfrac{a^{\sigma 0(\alpha,r,j)} - 1}{\ln a}$$

$$\sigma 07(\alpha,m,r,j) = m\left(\sqrt[m]{1+\sigma 0(\alpha,r,j)} - 1\right)$$

It should be taken into account that some of these multipliers may change when the parameter is changed $\alpha$. In addition, it should be noted that the convergence factors $\sigma 5(\alpha,a,r,j)$, $\sigma 6(\alpha,a,r,j)$, $\sigma 7(\alpha,m,r,j)$ depend on two parameters.

2. **The method of Fourier coefficients of class functions $W_v^r$.**

From the course of mathematical analysis it is known that the Fourier coefficients $a_j$, ($j=1,2,...$) of the even function $f(x) \in W_v^r$ are in descending order $j^{-(r+1)}$. From this it follows that in the role of multipliers of convergence with the order of descent $j^{-(r+1)}$ one can apply these Fourier coefficients. At the same time, the main difficulties lie in the construction of the function itself $f(x) \in W_v^r$.

Below we consider two types of such functions, which were considered in [12].

a) consider the function

$$T^*(r,\alpha,t) = \begin{cases} \left\{C_r\left[\sin\left[\dfrac{\pi}{2}\left(1-\left|\dfrac{t}{\alpha}\right|\right)\right]\right]\right\}^r, & \text{if } |t| \le \alpha < \pi; \\ 0, & \text{otherwise } t, \end{cases} \quad (r=1,2,...).$$

Parameter $C_r$ is determined from the condition of normalization of the function $T(r,\alpha,x)$,

$$C_r^{-1} = \int_{-\pi}^{\pi} T(r,\alpha,x)\,dx.$$

It is easy to make sure that $T^*(r,\alpha,t) \in W_v^r$. Calculating the coefficients of this function by cosines, we have:

for odd values $r = 2l+1$, ($l = 0,1,...$) we have

$$\mu_k(2l+1,\alpha) = \dfrac{(-1)^{l+1}\left[(2l+1)!!\right]^2 \pi^{2(l+1)} \cos(\alpha k)}{\left(4\alpha^2 k^2 - 1^2 \cdot \pi^2\right)\left(4\alpha^2 k^2 - 3^2 \cdot \pi^2\right)...\left(4\alpha^2 k^2 - (2l+1)^2 \cdot \pi^2\right)};$$

for even values $r = 2l$, ($l=1,2,...$) we have

$$\mu_k(2l,\alpha) = \dfrac{(-1)^l (l!)^2 \pi^{2l}}{\left(\alpha^2 k^2 - 1^2 \cdot \pi^2\right)\left(\alpha^2 k^2 - 2^2 \cdot \pi^2\right)...\left(\alpha^2 k^2 - l^2 \cdot \pi^2\right)} \dfrac{\sin \alpha k}{\alpha k}.$$

c) consider the function of the species

$$P(r,\alpha,t) = \begin{cases} C_r\left[1-\left(\dfrac{t}{\alpha}\right)^2\right]^r, & |t|\le\alpha \\ 0, & |t|\le\alpha \end{cases} \quad (r=1,2,\ldots)$$

where, as before, is a parameter $C_r$ is determined from the normalization conditions.

Calculating the Fourier coefficients of this function by cosines, we get for $r=1,2,3,4$

$$v_k(1,\alpha) = \frac{3(\sin(\alpha k) - \alpha k\cos(\alpha k))}{\pi\alpha^3 k^3}, \quad (k=1,2,\ldots).$$

$$v_k(2,\alpha) = \frac{-15\left((\alpha^2 k^2 - 3)\sin(\alpha k) + 3\alpha k\cos(\alpha k)\right)}{\pi\alpha^5 k^5};$$

$$v_k(3,\alpha) = \frac{105\left((15 - 6\alpha^2 k^2)\sin(\alpha k) - (15 - \alpha^2 k^2)\alpha k\cos(\alpha k)\right)}{\pi\alpha^7 k^7};$$

$$v_k(4,\alpha) = \frac{945\left[(105 - 45\alpha^2 k^2 + \alpha^4 k^4)\sin(\alpha k) - (105 - 10\alpha^2 k^2)\alpha k\cos(\alpha k)\right]}{\pi\alpha^9 k^9}.$$

Reasoning similarly, it is easy to obtain the Fourier coefficients for other values of the parameter $r$.

3. **Variable convergence factors of the type $\sigma(\alpha,r,j) = [\operatorname{sinc}(\alpha j)]^{1+r}$.**

These multipliers deserve special attention and have quite interesting properties, which we plan to consider in future works. Note that at $r=0$ multipliers $\sigma(\alpha,r,j)$ to some extent related to the Fourier transform of a rectangular pulse.

We present a general method of constructing functions of a system of simple fundamental trigonometric splines; these formulas have the form

$$st_k(I_1,I_2,\sigma,r,N,t) = \frac{1}{N}\left(1 + 2\sum_{j=1}^{\frac{N-1}{2}} \frac{c_j(I_1,I_2,\sigma,r,k,N,t)}{h_j(I_1,I_2,\sigma,r,N)}\right),$$

where

$$c_j(I_1,I_2,\sigma,r,k,t) = \sigma(r,j)\cos(j(t-(x_k^{(0)})^{1-I_2}(x_k^{(1)})^{I_2})) +$$

$$+ \sum_{m=1}^{\infty}(-1)^{m(r+1+I_1+I_2)}\left[\sigma(r,mN-j)\cos((mN-j)(t-(x_k^{(0)})^{1-I_2}(x_k^{(1)})^{I_2})) +\right.$$

$$+\sigma(r,mN+j)\cos((mN+j)(t-(x_k^{(0)})^{1-I_2}(x_k^{(1)})^{I_2}))\bigg],$$

$$h_j(I_1,I_2,\Gamma,\sigma,r,N) = \sigma(r,j) + \sum_{m=1}^{\infty}(-1)^{m(r+1+I_1+I_2)}\left[\sigma(r,mN-j) + \sigma(r,mN+j)\right].$$

Note that at $\alpha = \pi/N$ this method of constructing trigonometric fundamental simple splines with sign-changing factors $\sigma(\alpha,r,j)$ and in the given method of constructing trigonometric fundamental simple splines with constant coefficients $\sigma 0(\alpha,r,j)$ lead to the same results.

4. **The method of constructing sign-changing convergence factors, , is based on the equivalence of infinitesimals.**

As before, note that the coefficients of convergence are constant $\sigma(\alpha,r,j)$ at $j\to\infty$ are infinitely small quantities. Using the known equivalence relations of infinitesimals, it is easy to construct the following convergence factors:

$$\sigma 1(\alpha,r,j) = \sin(\sigma(\alpha,r,j));$$

$$\sigma 2(\alpha,r,j) = \arcsin(\sigma(\alpha,r,j));$$

$$\sigma 3(\alpha, r, j) = tg\left(\sigma(\alpha, r, j)\right);$$

$$\sigma 4(\alpha, r, j) = arctg\left(\sigma(\alpha, r, j)\right);$$

$$\sigma 5(\alpha, a, r, j) = \ln a \log_a\left(\sigma(\alpha, r, j)\right);$$

$$\sigma 6(\alpha, a, r, j) = \frac{a^{\sigma(\alpha, r, j)} - 1}{\ln a};$$

$$\sigma 7(\alpha, m, r, j) = m\left(\sqrt[m]{1 + \sigma(\alpha, r, j)} - 1\right).$$

5. **Method of composite factors of convergence.**

This method consists in the fact that when constructing the sum of the Fourier series, to which trigonometric splines are given, a composite convergence factor is used, which is given as follows

$$v(r, j) = \begin{cases} c_j, & 1 \le j \le N; \\ \varphi(r, j), & j > N, \end{cases}$$

where $c_j$, $(1 \le j \le N)$ - limited $\varphi(r, j) = O(j^{-(1+r)})$, a $1 < N < \infty$.

This method allows you to construct trigonometric splines of predefined shapes.

6. **Convergence multipliers obtained from the convergence multipliers of previous groups by algebraic operations.**

The multipliers of the previous groups can be multiplied, divided, taken modulo, etc. Of course, attention should be paid to the decreasing order of the multipliers obtained as a result of performing these operations, as well as to the parameters on which the original multipliers depended. So, for example, as a result of multiplying the factor of the 1st group by the factor of the 4th group, we get the factor

$$\sin\left(\frac{\alpha}{j^{1+r_1}}\right)\left(\frac{\sin(\beta j)}{\beta j}\right)^{1+r_2},$$

which already depends on two parameters $\alpha$ and $\beta$ in descending order $2 + r_1 + r_2$. It is also possible to consider the convergence factors of the type

$$(-1)^j \frac{\alpha}{j^{1+r}}, \frac{1+(-1)^j}{2} \frac{\alpha}{j^{1+r}},$$

etc.

**Conclusions.**
1. The connection between the class of polynomial interpolation simple splines and the class of trigonometric interpolation simple splines is established.
2. When studying trigonometric interpolation splines, it is necessary to consider stitching grids and interpolation grids, since various combinations of stitching grids and interpolation grids are possible when constructing these splines.
3. A general method of constructing trigonometric fundamental splines with vectors is given Γ and H and trigonometric fundamental simple splines with constant coefficients of convergence.
4. An example and graphs illustrating the introduced concepts of stitching grids and interpolation grids are given.
5. Methods of constructing one- and two-parameter sign-fixed and sign-changing convergence factors based on the equivalence of infinitesimals are proposed.
6. A method of constructing one-parameter convergence factors as coefficients of the Fourier series of class functions is proposed $W_v^r$.
7. The special role of sign-changing convergence factors is noted $\sigma(\alpha, r, j) = [\text{sinc}(\alpha j)]^{1+r}$, which have quite interesting properties.
8. A general method of constructing trigonometric fundamental simple splines with sign-changing convergence factors is given $\sigma(\alpha, r, j) = [\text{sinc}(\alpha j)]^{1+r}$.
9. The method of constructing composite convergence factors is considered, which allows constructing trigonometric splines of predefined shapes.

10. The method of constructing convergence factors by applying algebraic operations is considered.
11. Undoubtedly, the influence of the proposed convergence factors on the approximate properties of trigonometric splines requires further research.

## List of references